\documentclass[final,twoside,11pt]{entics}
\usepackage{enticsmacro}
\usepackage{graphicx}
\usepackage[all]{xy}
\def\uuar{\mathord{\mbox{\makebox[0pt][l]{\raisebox{.4ex}
{$\uparrow$}}$\uparrow$}}}

\def\da{\mathord{\downarrow}}

\def\R{{\mathbb R}}

\def\ur{(U,R)}
\def\r-{\underline{R}}
\def\R-{\overline{R}}

\def\dda{\mathord{\mbox{\makebox[0pt][l]{\raisebox{-.4ex}
{$\downarrow$}}$\downarrow$}}}

\def\da{\mathord{\downarrow}}

\def\R{{\mathbb R}}

\newtheorem{tm}{Theorem}[section]
\newtheorem{pn}[tm]{Proposition}
\newtheorem{lm}[tm]{Lemma}
\newtheorem{cy}[tm]{Corollary}

\newtheorem{dn}[tm]{Definition}
\newtheorem{rk}[tm]{Remark}
\newtheorem{ex}[tm]{Example}

\def\conf{ISDT 2022} 	
\volume{2}			
\def\volu{2}			
\def\papno{13}			
\def\mydoi{{\normalshape  \href{https://doi.org/10.46298/entics.10420}{10.46298/entics.10420}}}
\def\copynames{G. Wu, L. Xu} 
\def\runtitle{Representations of Domains via CF-approximation Spaces}

\def\lastname{Wu and Xu}

\begin{document}

 \begin{frontmatter}

\title{Representations of Domains via CF-approximation
Spaces} \author{Guojun Wu\thanksref{ALL}\thanksref{coemail}}     \author{Luoshan Xu\thanksref{ALL}\thanksref{myemail}}
  \address{College of Mathematical Science \\ Yangzhou University\\
   Yangzhou 225002, P. R. China}
\thanks[ALL]{Supported by National Natural Science Foundation of China (11671008).}
\thanks[coemail]{Email:
    \href{2386858989@qq.com} {\texttt{\normalshape
        2386858989@qq.com}}} \thanks[myemail]{Email:
    \href{mailto:mailto:luoshanxu@hotmail.com} {\texttt{\normalshape
        luoshanxu@hotmail.com}}}

\begin{abstract}
Representations of domains mean in a general way representing a domain as a suitable family endowed with set-inclusion order of some mathematical structures.  In this paper, representations of domains via CF-approximation
spaces are considered. Concepts of CF-approximation
spaces and CF-closed sets are introduced. It is proved that the family of CF-closed sets in a CF-approximation
space endowed with set-inclusion order is a continuous domain and that every continuous domain is  isomorphic to the  family of CF-closed sets of some CF-approximation
space endowed with set-inclusion order. The concept of CF-approximable relations is introduced using a categorical approach, which later facilitates the proof that the category of CF-approximation spaces and CF-approximable relations is equivalent to that of continuous domains and Scott continuous maps.
\end{abstract}
\begin{keyword} CF-approximation
space; CF-closed set; CF-approximable relation;  continuous domain; abstract base
\end{keyword}
\end{frontmatter}
\section{Introduction}
Domain theory   is one of the important research fields of theoretical computer science  \cite{Gierz}. In recent years of research in domain theory, there is a growing body of scholarly work towards synthesizing various mathematical fields such as ordered structures, topological spaces, formal contexts, rough sets, and various kinds of logic.
 One of such syntheses is to create representation for various kinds of domains using abstract bases \cite{AJung,wang2}, formal topologies \cite{XU-mao-form}, information systems \cite{Spreen-xu-inf,wang2}, formal contexts \cite{guo1}-\cite{wang}, and so on.  Amongst these, representation via abstract bases appears to be most natural due to its simplicity.

By representation of domains, we mean any general way by which one can characterize a domain using a suitable family of some mathematical structures ordered by the set-theoretic inclusion. With this understanding, clearly, every continuous domain can be represented by c-infs  \cite{Spreen-xu-inf}, abstract bases, formal topologies \cite{XU-mao-form}, etc.   Recently, Qingguo Li, et. al. in  \cite{wang} introduced attribute continuous formal contexts which are quadruples,  and showed that every continuous domain can be represented by  attribute continuous formal contexts.

 While representation using abstract bases appears to be the most natural and simple, its scope of study is unfortunately too narrow in that it is easy to miss out on something deeper. Noticing from rough set theory  \cite{jar2} that abstract bases are all special generalized approximation
spaces (GA-spaces, for short) \cite{Yang-Xu-Alg},  we consider   generalize  an abstract base to a CF-approximation
space  which is a GA-space  with  some coordinating family of finite sets. Since  the lower approximation operator $\underline{R}$ and the upper approximation operator $\overline{R}$ are mutually dual in a CF-approximation space, we mainly use the upper approximation operator  $\overline{R}$ and introduce CF-closed sets which are  generalizations  of round ideals in  abstract bases.
 With these concepts, representations of domains via CF-approximation
spaces are obtained. We will see that this approach of  representing  domains is more general than the approach of representing  domains by abstract bases.  We also introduce the concept of CF-approximable relations  using a categorical approach  and prove that the category of CF-approximation spaces and CF-approximable relations is equivalent to that of continuous domains and Scott continuous maps.  This work makes links naturally between domains and rough sets.

\section{Preliminaries}

We quickly recall some basic notions and results of domain theory. For a set $U$ and $X\subseteq U$,  we use $\mathcal {P} (U)$
to denote the power set of $U$, $\mathcal{P}_{fin}(U)$ to denote the family of all nonempty finite subsets of $U$  and  $X^c$ to denote the complement
of $X$ in $U$. The symbol $F\subseteq_{fin} X$ means $F$ is a finite subset of $X$. For  notions which we do not explicitly define herein, the reader may refer to  \cite{Gierz,Jean-2013}.

Let ($L$, $\leqslant$) be a poset. A \emph{principal ideal} (resp., \emph{principal filter}) of $L$
is a set of the form ${\da x}=\{y\in L\mid y\leqslant x\}$ (resp., ${\uparrow\! x}=\{y\in L\mid x\leqslant y\}$). For $A\subseteq L$, we write $\da A=\{y\in L\mid\exists \ x\in A, \,y\leqslant x\}$ and $\uparrow\! A=\{y\in L\mid\exists \ x\in A,\, x\leqslant
y\}$. A subset $A$ is a \emph{lower set} (resp., an \emph{upper set}) if $A=\da A$ (resp., $A=\uparrow\!A$). We say that $z$ is a \emph{lower bound} (resp., an \emph{upper bound}) of $A$ if $A\subseteq \uparrow\!z$ (resp., $A\subseteq\downarrow\!z$).  The supremum of $A$ is the least upper bound of $A$, denoted by $\bigvee A$ or $\sup A$. The infimum of $A$ is the greatest lower bound of $A$, denoted by $\bigwedge A$ or $\inf A$. A nonempty subset $D$ of $L$ is \emph{directed} if every  finite subset of $D$ has an upper bound in $D$. A  subset $C$ of $L$ is \emph{consistent} if $C$ has an upper bound in $L$.
 A poset $L$ is a \emph{directed complete partially ordered set} ({\em dcpo}, for short) if every directed subset of $L$ has a supremum. A \emph{semilattice} (resp., {\em sup-semilattice}) is a poset in which every pair of elements has an infimum (resp., a supremum). A \emph{complete lattice} is a poset in which every subset has a supremum (equivlently, has an infimum).
If any finite consistent subset $A$ of $L$ has a  supremum, then $L$ is called  a {\em cusl}. If any  consistent subset $B$ of  $L$ has a  supremum, then $L$ is called  a {\em bc-poset}.

Let $L$ be a semilattice and $K\subseteq L$. If for all $x, y \in K$ it holds that $x \wedge y \in K$, then $K$ is called a {\em subsemilattice} of $L$.

Recall that in a poset $P$, we say that $x$ \emph{way-below}
$y$, written $x\ll y$  if whenever $D$ is a directed set that has
a supremum with $\sup D\geqslant y$, then $x\leqslant d$ for some $d\in D$. If $x\ll x$, then $x$ is called a compact element of $P$, the set $\{x\in P\mid x\ll x\}$ is denoted by $K(P)$. 
The set $\{y\in P\mid  x\ll
y\}$ will be denoted by $\uuar x$ and $\{y\in P\mid  y\ll x\}$ denoted
by $\dda x$.  A poset $P$ is said to be
\emph{continuous} (resp., {\em algebraic}) if for all $ x\in P$, $\dda x$ is directed (resp., ${\downarrow} x\cap K(P)$ is directed) and  $x=\bigvee \dda x$ (resp., $x=\bigvee ({\downarrow} x\cap K(P))$). If a dcpo $P$ is  continuous (resp., algebraic), then $P$ is called a {\em continuous domain} (resp., {\em an algebraic domain}). If a continuous domain $P$ is a semilattice (resp., sup-semilattice, complete lattice), then $P$ is called a {\em continuous semilattice} (resp., {\em continuous sup-semilattice}, {\em continuous lattice}). If a bc-poset $P$ is also a continuous domain, then $P$ is called  a {\em bc-domain}.  If  an algebraic domain $L$ is a semilattice and $K(L)$ is  a subsemilattice of $L$, then $L$ is called {\em arithmetic semmilattice}.

Let $L$ and $P$ be dcpos, and $f: L \longrightarrow P$ a map. If for any directed subset $D\subseteq L$, $f(\bigvee D)=\bigvee f(D)$, then $f$ is called a {\em Scott continuous map}.

\begin{lm} {\em (\cite{Gierz})}\label{lm-way-rel} Let $P$ be a poset. Then for all $x, y, u, z\in P$,
\begin{enumerate}
\item[{\em (1)}] $x\ll y\Rightarrow x\leqslant y$; 
\item[{\em (2)}] $u\leqslant x\ll y\leqslant z\Rightarrow u\ll z$.
\end{enumerate}
\end{lm}
\begin{lm}\label{pn-int}
If $P$ is a continuous poset, then the way-below relation
$\ll$ has
the  interpolation property:
$$x\ll z\Rightarrow \exists y\in P \mathrm{\ such\ that\ } x\ll y\ll
   z.$$
\end{lm}
\begin{dn}\label{base}
 Let $P$ be a poset,
$B\!\subseteq\!P$. The set $B$ is called a basis for $P$ if for all
$ a\in\!P$, there is a directed set $D_a\subseteq\!B\cap\dda a$ such
that $\sup_PD_a=\!a$, where the subscripted $P$ indicates that the operation (in this case, the supremum) is taken  in poset $P$.
\end{dn}

 It is well known that a poset $P$ is continuous iff it has a
 basis and that $P$ is algebraic  iff $K(P)$ is a
 basis.

 A binary relation $R\subseteq U\times U$ on a set $U$ is called
\emph{transitive} if $xRy$ and $yRz$ implies $xRz$ for all
$x,y,z\in U$. A binary relation $R$ is said a {\em preorder} if it is  reflexive and transitive.
\begin{dn}  {\em (see  \cite{Gierz,Jean-2013})} Let $(U, \prec)$ be a set
 equipped with a binary relation. The
 binary relation $\prec$  is called fully transitive
if it is transitive
and satisfies the strong interpolation property:
\begin{center}
$\forall |F| <\infty, F\prec z\Rightarrow \exists y\prec z$ such
that $F\prec y$,
\end{center}
where $F\prec y$ means for all $ t\in F$, $t\prec y$.  If $(B,\prec
)$ is a set equipped with a binary relation
 which is fully transitive, then $(B,\prec )$ is called an abstract
basis.
 \end{dn}
\begin{dn}
{\em (\cite{Gierz,Jean-2013})} Let $(B, \prec )$ be an abstract basis. A
non-empty subset $I$ of $B$ is a round ideal if

$(1)$ $\forall y \in I$, $x\prec y \Rightarrow x\in I$;

$(2)$ $\forall x$, $y \in I$, $\exists z\in I$ such that
$x\prec z$ and $y\prec z$.

\noindent All the round ideals of $B$ in set-inclusion order  is
called the round ideal completion of $B$, denoted by $RI(B)$.
\end{dn}

Observe that if $B$ is a basis for a continuous domain $P$, then
$(B,\,\ll)$, the restriction of the way-below relation to $B$,
is an abstract basis. And it is known (see
in  \cite{Spreen-xu-inf}) that $P$ in this case is isomorphic to $RI(B)$.

\begin{pn}\label{cy-gbase-idl}
If $P$ is a continuous domain, then $(P,\ll)$ is an abstract basis and  $RI(P,\ll)\cong (P, \leqslant)$.
\end{pn}


Next, we introduce some terminologies imported from rough set theory. A set  $U$ with
a binary relation $R$ is called a \emph{generalized approximation space}
(\emph{GA-space}, for short). Let $\ur$ be a GA-space.
Define $R_s, R_p: U\rightarrow\mathcal {P}(U)$  such that for all $x\in U$,
$R_s(x)=\{y\in U \mid   xRy\}, R_p(x)=\{y\in U \mid  yRx\}.$

Lower and upper approximation operators are key notions in GA-spaces.
\begin{dn}{\em (cf.  \cite{w-zhu})} \label{dn-app}
 Let $(U,R)$ be a GA-space. For $A\subseteq U$,
 define

\begin{center}
$\underline{R} (A)=\{x\in U\mid \ R_s(x)\subseteq A\}$,\qquad
$\overline{R}(A)=\{x\in U\mid \ R_s(x)\cap A\neq \emptyset\}.$
\end{center}

\noindent The
operators $\underline{R},\overline{R}:\mathcal {P} (U)\rightarrow
\mathcal {P}(U)$  are respectively called the lower and upper
approximation operators in $(U,R)$, where $\mathcal {P} (U)$ is the
power set of $U$.
\end{dn}

\begin{lm}{\em (cf.  \cite{gl-liu})}\label{lm-ula}
Let $(U,R)$ be a GA-space. Then the lower and upper approximation
operators $\underline{R}$ and $\overline{R}$ have the following
properties.
\begin{enumerate}
\item[{\em(1)}] $\underline{R} (A^c)=(\overline{R}(A))^c$,
$\overline{R}(A^c)=(\underline{R}(A))^c$, where $A^c$ is the
complement of $A\subseteq U$.

\item[{\em(2)}] $\underline{R}(U)=U$, $\overline{R}(\emptyset)=\emptyset$.

\item[{\em(3)}] Let $\{A_{i}\mid i\in I\}\subseteq \mathcal {P}(U)$. Then
$\underline{R}(\bigcap_{i\in I}A_{i})=\bigcap_{i\in
I}\underline{R}(A_{i}),\ \overline{R}(\bigcup_{i\in
I}A_{i})=\bigcup_{i\in I}\overline{R}(A_{i}).$

\item[{\em(4)}] If $A\subseteq B\subseteq U$, then
$\underline{R}(A)\subseteq\underline{R}(B),\overline{R}(A)\subseteq\overline{R}(B)$.

\item[{\em(5)}] For all $x\in U, \ \overline{R}(\{x\})=R_p(x)$.
\end{enumerate}
\end{lm}
\begin{lm}\label{lm-ref}{\em \cite{yy-yao2}}
Let $(U,R)$ be a GA-space. Then $R$ is reflexive iff   for
all $X\subseteq U$,
 $X\subseteq \overline{R}(X)$.
\end{lm}

\begin{lm}\label{lm-tra}{\em \cite{yy-yao2}}
Let $(U,R)$ be a GA-space. Then $R$ is transitive iff for
all $X\subseteq U$,
 $\overline{R}(\overline{R}(X))\subseteq \overline{R}(X)$.
\end{lm}

\begin{dn} {\em(\cite{Yang-Xu-Top})}
Let $(U,R)$ be a GA-space and $A\subseteq U.$  The set $A$ is called
 {\em $R$-open} if $A\subseteq \underline{R}(A)$ and {\em $R$-closed} if $\overline{R}(B)\subseteq B$.
\end{dn}

For a  preorder $R$, the operator $\underline{R}$ is an interior operator of a topology, so we have
\begin{dn}{\em(\cite{Yang-Xu-Top})}  If $R$ is a preorder, then GA-space $(U, R)$ is called a topological GA-space.
\end{dn}



%
%
%

The next proposition shows that all $R$-open sets of $(U,R)$ form a topology on $U$.
\begin{pn}\label{pn-tp}
Let $(U,R)$ be a GA-space. Then $\tau_{R}=\{A\subseteq U\mid A\subseteq\underline{R}(A)\}$ is an Alexandrov topology on $U$. 
\end{pn}
\begin{proof}
  Obviously, $\emptyset, U\in\tau_R$. By Lemma \ref{lm-ula}(3), we have $\tau_R$ is closed under
arbitrary intersections.
Let $A_i\in\tau_R$, namely, $A_i\subseteq\underline{R}(A_i)$ $(i\in\!I)$.  Then  $\bigcup_{i\in\!I}A_i\subseteq\bigcup_{i\in\!I}(\underline{R}(A_i))$. It follows from  $\underline{R}(A_i)\subseteq \underline{R}(\bigcup_{i\in I}A_i)$ that $\bigcup_{i\in\!I}(\underline{R}(A_i))\subseteq\underline{R}(\bigcup_{i\in\!I}A_i)$. Then $\bigcup_{i\in\! I}A_i\subseteq\underline{R}(\bigcup_{i\in\!I}A_i)$. So $\bigcup_{i\in\!I}A_i\in\tau_R$, namely, $\tau_R$ is closed under arbitrary unions.
This shows that $\tau_{R}$ is an Alexandrov topology.
\end{proof}

The above topology $\tau_{R}$ is called a \emph{topology induced by relation} $R$. Obviously, all the R-closed sets of $(U,R)$ are precisely all the closed sets of $\tau_{R}$ .
%
%
%
%

\section{CF-approximation
Spaces and CF-closed Sets}

 For an abstract base $(B, \prec)$, we naturally have the  triple $(B, \prec, \{\{b\}\mid b\in B\})$ and that the family $\{{\downarrow^{\prec} b}\mid b\in B\}$ is a base of the continuous domain $RI(B)$, where $\downarrow^{\prec} b=\{c\in B\mid c\prec b\}$.  We generalize an abstract base to a GA-space with consistent family of finite subsets (CF-approximation space, for short) by changing $(B, \prec)$ to a GA-space $(U, R)$ with $R$ being transitive and changing the family $\{\{b\}\mid b\in B\}$ to a  suitable family $\mathcal{F}$ of some finite subsets of $U$. We hope that  the family $\mathcal{F}$ can also induce a base of a continuous domain.

\begin{dn}\label{dn-cf-ga}Let $(U, R)$ be a GA-space, $R$  a transitive relation and $\mathcal{F}\subseteq \mathcal{P}_{fin}(U)\cup\{\emptyset\}$. If for all $ F\in \mathcal{F}$,  whenever $K\subseteq_{fin} \overline{R}(F)$, there always exists $G\in \mathcal{F}$ such that $K\subseteq \overline{R}(G)$ and $G\subseteq \overline{R}(F)$, then $(U, R, \mathcal{F})$ is called a generalized approximation space with consistent family of finite subsets, or a CF-approximation space, for short.
\end{dn}
\begin{lm}\label{lm2-tr-up} Let $(U, R)$ be a GA-space, $A, B\subseteq U$. If $R$ is a transitive relation, then $\overline{R}(B)\subseteq\overline{R}(A)$
when $B\subseteq\overline{R}(A)$.
\end{lm}
\begin{proof} It follows from Lemma \ref{lm-ula}(4) and \ref{lm-tra}. \end{proof}

\begin{dn}\label{dn-cf-cl}Let $(U, R, \mathcal{F})$ be a CF-approximation space, $E\subseteq U$. If for all $ K\subseteq_{fin}E$, there always exists $F\in \mathcal{F}$ such that $K\subseteq \overline{R}(F)\subseteq E$ and $F\subseteq E$, then $E$ is called a CF-closed set of $(U, R, \mathcal{F})$. The collection of all CF-closed sets of $(U, R, \mathcal{F})$ is denoted by $\mathfrak{C}(U, R, \mathcal{F})$.
\end{dn}
\begin{rk} $(1)$ If $\emptyset\in \mathfrak{C}(U, R, \mathcal{F})$, then $\emptyset\in\mathcal{F}$ by $\overline{R}(\emptyset)=\emptyset$.

$(2)$ For CF-approximation space $(U, R, \mathcal{F})$, if $\mathcal{F}=\{\{x\}\mid x\in U\}$, then   $(U, R)$ is an abstract base by Lemma \ref{lm-ula}(5), and all the CF-closed sets of $(U, R, \mathcal{F})$ are precisely all the round ideals of  $(U, R)$.
\end{rk}

The following example shows that $(U, R)$ is  not necessarily an abstract base when $(U, R, \mathcal{F})$ is  a CF-approximation space, showing that CF-approximation spaces is a generalization of abstract bases.
\begin{ex} Let $U=\mathbb{N}$, $R=\{(1, 1), (1, 2), (1, 3), (1, 4), (2, 3), (2, 4), (4, 3)\}\cup\{(i, i)\mid i\geqslant 5\}$, $\mathcal{F}=\{\{1\}, \{1, 2\}, \emptyset\}\cup\{\{i\}\mid i\geqslant 5\}$. It is easy to check that $(U, R, \mathcal{F})$ is  a CF-approximation space and $\mathfrak{C}(U, R, \mathcal{F})=\{\emptyset, \{1\}\}\cup\{\{i\}\mid i\geqslant 5\}$ is a continuous domain. Notice that $(1, 4), (2, 4)\in R$, but there is no $t\in U$ such that $(1, t), (2, t), (t, 4)\in R$. So $(U, R)$ is not an abstract base.
\end{ex}

\begin{pn}\label{pn1-cf-clo}Let $(U, R, \mathcal{F})$ be a CF-approximation space. If $E\in \mathfrak{C}(U, R, \mathcal{F})$, then $E$ is an $R$-closed set.
\end{pn}
\begin{proof}
If $x\in \overline{R}(E)$, then $R_{s}(x)\cap E\neq\emptyset$. So there is $y\in U$ such that $xRy$ and $y\in E$. By Definition \ref{dn-cf-cl}, there exists $F\in \mathcal{F}$, such that $y\in \overline{R}(F)\subseteq E$ and $F\subseteq E$. By the transitivity of $R$, we know that $\overline{R}(\{y\})\subseteq \overline{R}(\overline{R}(F))\subseteq \overline{R}(F)\subseteq E$. Thus $\overline{R}(\{y\})\subseteq E$. It is clear  that $x\in \overline{R}(\{y\})\subseteq E$ because of $xRy$.  By the arbitrariness of $x\in \overline{R}(E)$, we know that $\overline{R}(E)\subseteq E$.  This shows that $E$ is an $R$-closed set. \end{proof}

\begin{pn}\label{pn2-cf-clo} Let $(U, R, \mathcal{F})$ be a CF-approximation space,  then the following statements hold:
\begin{enumerate}
\item[$(1)$] For any $F\in \mathcal{F}$, $\overline{R}(F)\in \mathfrak{C}(U, R, \mathcal{F})$;

\item[$(2)$] If $E\in \mathfrak{C}(U, R, \mathcal{F})$, $A\subseteq E$, then $\overline{R}(A)\subseteq E$;

\item[$(3)$] If $\{E_{i}\}_{i\in I}\subseteq\mathfrak{C}(U, R, \mathcal{F})$ is a directed family, then $\bigcup_{i\in I}E_{i}\in \mathfrak{C}(U, R, \mathcal{F})$.
\end{enumerate}
\end{pn}
\begin{proof}
(1)    Follows directly by Definition \ref{dn-cf-ga} and the transitivity of $R$.

(2) Follows from $A\subseteq E$ and Lemma \ref{lm-ula}(4) that $\overline{R}(A)\subseteq \overline{R}(E)$. By Proposition \ref{pn1-cf-clo}, we know that  $\overline{R}(E)\subseteq E$. Thus $\overline{R}(A)\subseteq E$.

(3)  Follows directly from Definition \ref{dn-cf-cl}. \end{proof}

The proposition above shows that $(\mathfrak{C}(U, R, \mathcal{F}), \subseteq)$ is a dcpo. The following  proposition gives  equivalent characterizations of  CF-closed sets.
\begin{pn}\label{pn3-cf-cl}
The following statements are equivalent for a CF-approximation space $(U, R, \mathcal{F})$:
\begin{enumerate}
\item[$(1)$] $E\in \mathfrak{C}(U, R, \mathcal{F})$;

\item[$(2)$] The family $\mathcal{A}=\{\overline{R}(F)\mid F\in \mathcal{F},  F\subseteq E\}$ is directed and $E=\bigcup\mathcal{A}$;

\item[$(3)$] There exists a family $\{F_{i}\}_{i\in I}\subseteq \mathcal{F}$ such that $\{\overline{R}(F_{i})\}_{i\in I}$ is directed, and $E=\bigcup_{i\in I}\overline{R}(F_{i})$;

\item[$(4)$]  There always  exists $F\in \mathcal{F}$ such that $K\subseteq \overline{R}(F)\subseteq E$ whenever $ K\subseteq_{fin}E$.
\end{enumerate}
\end{pn}
\begin{proof} If $E=\emptyset\in  \mathfrak{C}(U, R, \mathcal{F})$, the proposition holds obviously. Let $E\neq\emptyset$.

$(1)\Rightarrow (2)$ By  Definition \ref{dn-cf-cl}, we know that $\mathcal{A}$ is not empty. Let $X_{1}, X_{2}\in \mathcal{A}$, then there exist $F_{1}, F_{2}\in \mathcal{F}$ and $F_{1}, F_{2}\subseteq E$, such that $X_{1}=\overline{R}({F_{1}})$, $X_{2}=\overline{R}({F_{2}})$. By $F_{1}\cup F_{2}\subseteq_{fin} E$ and Definition \ref{dn-cf-cl} we know that there exists $F_{3}\in \mathcal{F}$, such that $F_{1}\cup F_{2}\subseteq\overline{R}(F_{3})$ and $F_{3} \subseteq E$.
 By the transitivity of $R$ and Lemma \ref{lm2-tr-up} we know that $\overline{R}(F_{1})\subseteq\overline{R}(F_{3})$, $\overline{R}(F_{2})\subseteq\overline{R}(F_{3})$. This shows that $\mathcal{A}$ is directed. Next we prove $E=\bigcup\mathcal{A}$. By Proposition \ref{pn2-cf-clo}(2) we know that $\bigcup\mathcal{A}\subseteq E$ holds. Conversely, if $x\in E$, then by Definition \ref{dn-cf-cl}, there is $F\in \mathcal{F}$ such that $x\in \overline{R}(F)\subseteq E$ and $F\subseteq E$. So $x\in \bigcup\mathcal{A}$. By the arbitrariness of $x\in E$ we know that $E\subseteq\bigcup\mathcal{A}$. Thus $E=\bigcup\mathcal{A}$.

$(2)\Rightarrow (3)$ Trivial.

$(3)\Rightarrow (4)$ Follows directly from  the finiteness of $ K$ and the directedness of  $\{\overline{R}(F_{i})\}_{i\in I}$.

$(4)\Rightarrow (1)$   If $ K\subseteq_{fin}E$, then there exists $F\in \mathcal{F}$ such that $K\subseteq \overline{R}(F)\subseteq E$.  By  Definition \ref{dn-cf-ga}, there exists $G\in \mathcal{F}$ such that $K\subseteq \overline{R}(G)$ and $G\subseteq \overline{R}(F)$. By Lemma \ref{lm2-tr-up} we know that  $\overline{R}(G)\subseteq\overline{R}(F)\subseteq E$. Thus $K\subseteq \overline{R}(G)\subseteq E$. Noticing that $G\subseteq \overline{R}(F)\subseteq E$,  by Definition \ref{dn-cf-cl} we obtain that $E\in \mathfrak{C}(U, R, \mathcal{F})$.
\end{proof}

The following theorem characterizes the way-below relation $\ll$ in  dcpo $(\mathfrak{C}(U, R, \mathcal{F}), \subseteq)$.

\begin{tm}\label{tm-cfga-way} Let $(U, R, \mathcal{F})$ be a CF-approximation space, $E_{1}, E_{2}\in \mathfrak{C}(U, R, \mathcal{F})$. Then $ E_{1}\ll E_{2}$ if and only if there exists $F\in \mathcal{F}$ such that $ E_{1}\subseteq \overline{R}(F)$ and  $F\subseteq E_{2}$.
\end{tm}
\begin{proof}
$\Rightarrow$:  It follows from $E_{2}\in \mathfrak{C}(U, R, \mathcal{F})$ and  Proposition \ref{pn3-cf-cl}(2)  that $E_{2}=\bigcup\{\overline{R}(F)\mid F\in \mathcal{F}, F\subseteq E_{2}\}$  and that $\{\overline{R}(F)\mid F\in \mathcal{F}, F\subseteq E_{2}\}$ is directed. If $ E_{1}\ll E_{2}$, then there exists $F\in \mathcal{F}$ such that $F\subseteq E_{2}$, $ E_{1}\subseteq \overline{R}(F)$.

$\Leftarrow$: For any directed family $\{C_{i}\}_{i\in I}\subseteq\mathfrak{C}(U, R, \mathcal{F})$, if $ E_{2}\subseteq \bigvee_{i\in I}C_{i}= \bigcup_{i\in I}C_{i}$, then by $F\subseteq E_{2}$ and the finiteness of $F$ we know that there exists $i_{0}\in I$ such that $F\subseteq C_{i_{0}}$. By Proposition \ref{pn2-cf-clo}(2)  and $E_{1}\subseteq \overline{R}(F)$, we know that $E_{1}\subseteq \overline{R}(F)\subseteq C_{i_{0}}$, showing that $E_{1}\ll E_{2}$.
\end{proof}

\begin{cy}\label{cy-cfga-way} Let $(U, R, \mathcal{F})$ be a CF-approximation space, $E\in \mathfrak{C}(U, R, \mathcal{F})$, $F\in \mathcal{F}$. The following statements hold:
\begin{enumerate}
\item[$(1)$] If $F\subseteq E$, then $\overline{R}(F)\ll E$;

\item[$(2)$] $\overline{R}(F)\ll \overline{R}(F)$ if and only if there exists $G\in \mathcal{F}$, such that $G\subseteq\overline{R}(G)=\overline{R}(F)$.
\end{enumerate}
\end{cy}
\begin{proof}
Follows directly from Theorem \ref{tm-cfga-way}.\end{proof}

\begin{tm}\label{tm-CF-CDO}Let $(U, R, \mathcal{F})$ be a CF-approximation space. Then $(\mathfrak{C}(U, R, \mathcal{F}), \subseteq)$ is a continuous domain.
\end{tm}
\begin{proof}
By Proposition \ref{pn2-cf-clo} we see that $(\mathfrak{C}(U, R, \mathcal{F}), \subseteq)$ is a dcpo. Set $\mathcal{B}=\{\overline{R}(F)\mid F\in \mathcal{F}\}$. Then $\mathcal{B}$ is a base of $(\mathfrak{C}(U, R, \mathcal{F}), \subseteq)$ by Proposition
\ref{pn3-cf-cl}(2) and Corollary \ref{cy-cfga-way}(1). Thus $(\mathfrak{C}(U, R, \mathcal{F}), \subseteq)$ is a continuous domain. \end{proof}


 The following theorem shows that any continuous domain $(L, \leqslant)$ can induce a CF-approximation space.

\begin{tm}\label{tm-CDO-CFGA} Let $(L, \leqslant)$ be a continuous domain, $R_{L}$  the way-below relation ``$\ll$" of $(L, \leqslant)$; $\mathcal{F}_{L}=\{F\subseteq_{fin} L\mid F\  \mbox{has a top element}\}$. For any $F\in \mathcal{F}_{L}$, let $c_{F}$ be the top element of $F$. Then $(L, R_{L},  \mathcal{F}_{L})$ is a CF-approximation space.
\end{tm}
\begin{proof}   By Lemma \ref{lm-way-rel}, we know that $R_{L}=\ll$ is transitive. For any $F\in \mathcal{F}_{L}$, by Lemma \ref{lm-ula}(5), we have that $\overline{R}_{L}(F)=\dda c_{F}$.  For $K\subseteq_{fin}\overline{R}_{L}(F)=\dda c_{F}$, by that $L$ is a continuous domain, we know that $\dda c_{F}$ is directed. Then there exists $x\in \dda c_{F}$ such that $K\subseteq {\downarrow} x$. It follows from  $x\ll c_{F}$ and Lemma \ref{pn-int} that there is $y\in L$ such that $x\ll y\ll c_{F}$. Thus $K\subseteq \dda y$. Set $G=\{y\}\in \mathcal{F}_{L}$. By that $K\subseteq \overline{R}_{L}(G)=\dda y$ and $G\subseteq\overline{R}_{L}(F)=\dda c_{F}$, we have that $(L, R_{L},  \mathcal{F}_{L})$ is a CF-approximation space. \end{proof}
%
\begin{tm}\label{tm-CDO-CFCLO}  Let $(L, \leqslant)$ be a continuous domain, $(L, R_{L},  \mathcal{F}_{L})$ the one constructed in Theorem \ref{tm-CDO-CFGA}. Then $\mathfrak{C}(L, R_{L},  \mathcal{F}_{L})=\{\dda x\mid x\in L\}$.
\end{tm}
\begin{proof} By the proof of Theorem \ref{tm-CDO-CFGA} and Proposition \ref{pn2-cf-clo}(1), we know that $\{\dda x\mid x\in L\}\subseteq \mathfrak{C}(L, R_{L},  \mathcal{F}_{L})$. Conversely, let $E\in \mathfrak{C}(L, R_{L},  \mathcal{F}_{L})$. Then by Proposition \ref{pn3-cf-cl}, there is  a directed set $D\subseteq L$ such that $E=\bigcup\{\dda d\mid d\in D\}$.  Next we prove $E=\dda \bigvee D$. Obviously, $E\subseteq \dda \bigvee D$. Conversely, if $x\in \dda \bigvee D$, then by Lemma \ref{pn-int},  there is  $y\in L$ such that $x\ll y\ll\bigvee D$. So there is $d\in D$ such that $x\ll y\leqslant d$. Thus $x\in \dda d\subseteq E$, and  $E=\dda \bigvee D$. This shows that $\mathfrak{C}(L, R_{L},  \mathcal{F}_{L})\subseteq\{\dda x\mid x\in L\}$ and $\mathfrak{C}(L, R_{L},  \mathcal{F}_{L})=\{\dda x\mid x\in L\}$. \end{proof}

\begin{tm}\label{tm-rep} {\em (Representation Theorem)} Let $(L, \leqslant)$ be a poset. Then $L$ is a continuous domain iff there exists a CF-approximation space $(U, R, \mathcal{F})$ such that $(L, \leqslant)\cong (\mathfrak{C}(U, R, \mathcal{F}), \subseteq))$.
\end{tm}
\begin{proof}
$\Leftarrow$: Follows directly by Theorem \ref{tm-CF-CDO}.

$\Rightarrow$: If $L$ is a continuous domain, then by Theorem \ref{tm-CDO-CFGA} we know that $(L, R_{L},  \mathcal{F}_{L})$ is a CF-approximation space. Define a map $f: L\to \mathfrak{C}(L, R_{L},  \mathcal{F}_{L})$ such that for all $x\in L$, $f(x)=\dda x$. Then it follows from  Theorem \ref{tm-CDO-CFCLO} and the continuity of $L$ that $f$ is an order
isomorphism.
\end{proof}

\section {Representations of some special domains}
In this section, we add some conditions to  CF-approximation spaces, and then  discuss  representations of some special types of continuous domains.
\begin{lm}\label{lm-bas-dom} Let $(L, \leqslant)$ be a continuous domain and $B\subseteq L$  a base. If $(B, \leqslant)$ is a semilattice (resp., sup-semilattice,  poset with bottom element,  poset with top element,  sup-semilattice with bottom element, cusl), then $(L, \leqslant)$ is a continuous semilattice (resp., continuous sup-semilattic,  continuous domain with bottom element,  continuous domain with  top element, continuous lattice, bc-domain).
\end{lm}
\begin{proof} (1) Let $(B, \leqslant)$ be a  semilattice. For any $x, y\in L$, set $D=\{a\wedge_{B}b\mid a\in \dda x\cap B, b\in \dda y\cap B\}$.  It is easy to show  that $\dda x\cap B$ and $\dda y\cap B$ are directed. Therefore $D$ is directed and $\bigvee D$ exists. It is clear that $\bigvee D\leqslant \bigvee  (\dda x\cap B)=x$, $\bigvee D\leqslant \bigvee(\dda y\cap B ) =y$. If $z\leqslant x, y $,  then for any $t\in \dda z\cap B$, we have $t\ll x=\bigvee  (\dda x\cap B)$, $t\ll y=\bigvee  (\dda y\cap B)$. Therefore there exist $t_{1}\in \dda x\cap B$, $t_{2}\in \dda y\cap B$ such that $t\leqslant t_{1}, t_{2}$. Thus $t\leqslant t_{1}\wedge_{B} t_{2}$. Noticing that $t_{1}\wedge_{B} t_{2}\in D$ and the arbitrariness of $t\in \dda z\cap B$, we have that  $z=\bigvee(\dda z\cap B)\leqslant \bigvee D$. This shows that $\bigvee D$ is a greatest lower bound of $x, y$, namely, $x\wedge y=\bigvee D$. Thus $L$ is a  continuous semilattice.

(2) Let $(B, \leqslant)$ be a sup-semilattice. For any $x, y\in L$, set $D=\{a\vee_{B}b\mid a\in \dda x\cap B, b\in \dda y\cap B\}$. Clearly, $D$ is directed and $\bigvee D$ exists. It is obvious that $x, y\leqslant \bigvee D$. Let $x, y\leqslant z$. Then for all $ a\in \dda x\cap B$ and $b\in \dda y\cap B$, we have that $a\ll z$, $b\ll z$.  By the directedness of $\dda z\cap B$, there exists $t\in  \dda z\cap B$ such that $a, b\leqslant t$. Therefore $a\vee_{B}b\leqslant  t$. Noticing that $a\vee_{B}b\in D$, we have $\bigvee D\leqslant \bigvee ( \dda z\cap B)=z$. This shows that $\bigvee D$ is a least upper bound of $x, y$, namely, $x\vee y=\bigvee D$. Thus $L$ is a continuous sup-semilattice.

(3)/(4)  If $\perp$/$\top$ is a bottom/top element of $(B, \leqslant)$, then $\perp$/$\top$ is also a bottom/top element of $L$.

(5) Let $(B, \leqslant)$ be a sup-semilattice with bottom element $\perp$. Then by (2) and (3), we know that $L$ is a sup-semilattice with bottom element. Since $L$ is a dcpo, $L$ is a complete lattice. Thus $L$ is a continuous lattice.

(6) Let $(B, \leqslant)$ be a cusl. For any $x, y, z \in L$ which satisfy $x, y\leqslant z$, we show that  for all $a\in \dda x\cap B$ and $ b\in \dda y\cap B$, $a\vee_{B}b$ exists. By that $x, y\leqslant z$, we have $a\ll z, b\ll z$. Since $B$ is a base, there exists $c\in\dda z\cap B$ such that $a, b\leqslant c$. That $a\vee_{B}b$ exists  by that $(B, \leqslant)$ is a cusl. Similar to the proof of (2), we have that $L$ is a cusl. As $L$ is a continuous domain, we see that $L$ is a bc-domain. \end{proof}

\begin{tm} Let $(U, R, \mathcal{F})$ be a CF-approximation space. If  $(\{\overline{R}(F)\mid F\in \mathcal{F}\}, \subseteq)$ is a semilattice (resp., sup-semilattice,  poset with bottom element,  poset with top element,  sup-semilattice with  bottom element, cusl), then $\mathfrak{C}(U, R, \mathcal{F})$ is a continuous  semilattice (resp., continuous sup-semilattic,  continuous domain with bottom element,  continuous domain with  top element, continuous lattice, bc-domain).  Conversely, any continuous semilattice (resp., continuous sup-semilattic,  continuous domain with  bottom element,  continuous domain with  top element, continuous lattice, bc-domain) $L$  is  isomorphic to $(\mathfrak{C}(L, R_{L},  \mathcal{F}_{L}), \subseteq)$ of  corresponding CF-approximation spaces, respectively.
\end{tm}
\begin{proof}
The first half of the theorem follows directly from Theorem \ref{tm-CF-CDO} and Lemma \ref{lm-bas-dom}.

For the second half, let $(L, R_{L},  \mathcal{F}_{L})$ be the one in  Theorem \ref{tm-CDO-CFGA}. Define a map $f: L\longrightarrow \{\overline{R}(F)\mid F\in \mathcal{F}_{L}\}$ such that for all $ x\in L$, $f(x)=\dda x$. Since $L$ is continuous, $f$ is an order isomorphism. Thus $\{\overline{R}(F)\mid F\in \mathcal{F}_{L}\}$ is a  semilattice (resp., sup-semilattice,  poset with bottom element,  poset with top element,  sup-semmilattice with  bottom element, cusl) whenever $L$ is a continuous  semilattice (resp., continuous sup-semilattic,  continuous domain with  bottom element,  continuous domain with  top element, continuous lattice, bc-domain). By Theorem \ref{tm-CDO-CFCLO}, $L$  is  isomorphic to $(\mathfrak{C}(L, R_{L},  \mathcal{F}_{L}), \subseteq)$.
\end{proof}

Next, we consider algebraic cases.

\begin{dn}\label{dn-cf-tga}  Let $(U, R)$ be a GA-space, $\mathcal{F}\subseteq \mathcal{P}_{fin}(U)\cup\{\emptyset\}$.  If $R$ is a preorder, then $(U, R, \mathcal{F})$ is called a topological CF-approximation space.
\end{dn}
\begin{rk} A topological CF-approximation space must be a CF-approximation space. In fact,   for all $ F\in \mathcal{F}$ and $K\subseteq_{fin} \overline{R}(F)$, taking $G=F$, then by Lemma \ref{lm-ref} we have  $K\subseteq \overline{R}(G)$ and $G\subseteq \overline{R}(F)$. By  Definition \ref{dn-cf-ga}, $(U, R, \mathcal{F})$ is a  CF-approximation space.
\end{rk}
\begin{pn}\label{pn-cfga-way} Let $(U, R, \mathcal{F})$ be a topological CF-approximation space, $E_{1}, E_{2}\in \mathfrak{C}(U, R, \mathcal{F})$. Then $ E_{1}\ll E_{2}$ iff there exists $F\in \mathcal{F}$  such that $ E_{1}\subseteq \overline{R}(F)\subseteq E_{2}$. Thus $K((\mathfrak{C}(U, R, \mathcal{F}), \subseteq))=(\{\overline{R}(F)\mid F\in \mathcal{F}\}, \subseteq)$.
\end{pn}
\begin{proof} It follows directly from Lemma \ref{lm-ref} and Theorem \ref{tm-cfga-way}. \end{proof}

\begin{tm}\label{tm-tpf-agl} Let $(U, R, \mathcal{F})$ be a topological CF-approximation space. Then $(\mathfrak{C}(U, R, \mathcal{F}), \subseteq)$ is an algebraic domain. Conversely,  any algebraic domain  can be represented by some topological CF-approximation space.
\end{tm}
\begin{proof} 
The first half of the theorem follows  from  Proposition \ref{pn-cfga-way} and   Theorem \ref{pn3-cf-cl} that $(\mathfrak{C}(U, R, \mathcal{F}), \subseteq)$ is an algebraic domain.

For the second half,  let $(L, \leqslant)$ be an algebraic domain. Set $(K(L),  R_{K(L)}, \mathcal{F}_{K(L)})$, where $\mathcal{F}_{K(L)}=\{F\subseteq_{fin} K(L)\mid F\ \mbox{has top element}\}$,  $R_{K(L)}=\leqslant$ is a partial order. Thus $(K(L),  R_{K(L)}, \mathcal{F}_{K(L)})$ is a topological CF-approximation space. For any $F\in \mathcal{F}_{K(L)}$, let $c_{F}$ be the top element of $F$.  By Lemma \ref{lm-ula}(5), we know that  for all $ F\in  \mathcal{F}_{K(L)}$, $\overline{R_{K(L)}}(F)={\downarrow} c_{F}\cap K(L)$. Similar to the proof of  Theorem \ref{tm-CDO-CFCLO}, we have that $\mathfrak{C}(K(L),  R_{K(L)}, \mathcal{F}_{K(L)})=\{{\downarrow}x\cap K(L)\mid x\in L\}.$  Since  $L$ is an algebraic domain, we know that  $(\{{\downarrow}x\cap K(L)\mid x\in L\}, \subseteq)\cong (L, \leqslant)$. The proof is finished. \end{proof}

\begin{lm}\label{lm-asl} Let $(L, \leqslant)$  be an algebraic domain. If $(K(L), \leqslant)$ is a semilattice, then $L$ is an arithmetic semilattice.
\end{lm}
\begin{proof} 
For any $x, y\in L$, let $D=\{a\wedge_{K(L)}b\mid a\in \downarrow x\cap K(L), b\in \downarrow y\cap K(L)\}$.  By Lemma \ref{lm-bas-dom}, we see that $x\wedge y=\bigvee D$ and $L$ is a semilattice. Next we prove $(K(L), \leqslant)$ is a subsemilattice of $L$. If $x, y\in K(L)$, then $x\wedge_{K(L)}y\in D$ and  $x\wedge_{K(L)}y$ is the top element of $D$. So $x\wedge_{K(L)}y=\bigvee D=x\wedge y$. Hence $x\wedge_{K(L)}y=x\wedge y$. This shows that $(K(L), \leqslant)$ is  a subsemilattice of $L$,  and $L$ is an arithmetic semilattice. \end{proof}

\begin{tm}\label{tm-CMF-CSL} Let $(U, R, \mathcal{F})$ be a topological CF-approximation space.  If poset $\{\overline{R}(F)\mid F\in \mathcal{F}\}, \subseteq)$ is a semilattice, then $(\mathfrak{C}(U, R, \mathcal{F}), \subseteq)$ is an arithmetic semilattice.  Conversely, any arithmetic semilattice can be represented in this way.
\end{tm}
\begin{proof}
The first half of the theorem follows  directly from  Theorem \ref{tm-tpf-agl} and Lemma \ref{lm-asl}.

For the second half,  let $L$ be an arithmetic semilattice and $(K(L),  R_{K(L)}, \mathcal{F}_{K(L)})$ be the one  in   Theorem \ref{tm-tpf-agl}. Therefore $(K(L),  R_{K(L)}, \mathcal{F}_{K(L)})$ is a topological CF-approximation space. For $F_{1}, F_{2}\in \mathcal{F}_{K(L)}$, then $\overline{R_{K(L)}}(F_{1})={\downarrow} c_{F_{1}}\cap K(L)$, $\overline{R_{K(L)}}(F_{2})={\downarrow} c_{F_{2}}\cap K(L)$. Since $L$ is an arithmetic semilattice, we know that $c_{F_{1}}\wedge c_{F_{2}}=c\in K(L)$.
So $\overline{R_{K(L)}}(F_{1})\cap\overline{R_{K(L)}}(F_{2})={\downarrow} c\cap K(L)$. This shows that there exists $\{c\}\in  \mathcal{F}_{L}$ such that $\overline{R_{K(L)}}(F_{1})\wedge\overline{R_{K(L)}}(F_{2})=\overline{R_{K(L)}}(\{c\})$. This shows that $\{\overline{R_{K(L)}}(F)\mid F\in \mathcal{F}\}$ is a semilttice. By Theorem
\ref{tm-tpf-agl}, we see that the second half of the theorem holds. \end{proof}

\section{CF-approximable Relations and Equivalence of Categories}
In this section, we define and study CF-approximable relations between  CF-approximation spaces and prove that the category of CF-approximation spaces and CF-approximable relations is equivalent to the category of continuous domains and Scott continuous maps.
\begin{dn}\label{dn-CF-app}  Let $(U_{1}, R_{1}, \mathcal{F}_{1})$, $(U_{2}, R_{2}, \mathcal{F}_{2})$ be CF-approximation spaces, and $\mathrel{\Theta}\subseteq\mathcal{F}_{1}\times\mathcal{F}_{2}$ a binary relation.  If
\begin{enumerate}
\item[$(1)$]  for all $F\in \mathcal{F}_1$, there is $G\in \mathcal{F}_{2}$ such that $F\mathrel{\Theta} G$;

\item[$(2)$]  for all $F, F^{\prime}\in \mathcal{F}_{1}$, $G\in  \mathcal{F}_{2}$, if $F\subseteq \overline{R_{1}}(F^{\prime})$, $F\mathrel{\Theta} G$, then $F^{\prime}\mathrel{\Theta} G$;

\item[$(3)$]  for all $F\in \mathcal{F}_{1}$, $G, G^{\prime}\in  \mathcal{F}_{2}$, if $F\mathrel{\Theta} G$, $G^{\prime}\subseteq \overline{R_{2}}(G)$,  then $F\mathrel{\Theta} G^{\prime}$;

\item[$(4)$] for all $F\in \mathcal{F}_{1}$,  $G\in \mathcal{F}_{2}$, if $F\mathrel{\Theta} G$, then there are $F^{\prime}\in \mathcal{F}_{1}$,  $G^{\prime}\in \mathcal{F}_{2}$ such that  $F^{\prime}\subseteq \overline{R_{1}}(F)$, $G\subseteq \overline{R_{2}}(G^{\prime})$ and $F^{\prime}\mathrel{\Theta} G^{\prime}$; and

\item[$(5)$]  for all $F\in \mathcal{F}_{1}$, $G_{1}, G_{2}\in \mathcal{F}_{2}$, if $F\mathrel{\Theta} G_{1}$ and $F\mathrel{\Theta} G_{2}$, then there is  $G_{3}\in \mathcal{F}_{2}$ such that  $G_{1}\cup G_{2}\subseteq\overline{R_{2}}(G_{3})$ and $F\mathrel{\Theta} G_{3}$,
\end{enumerate}
\noindent then $\mathrel{\Theta}$ is called a CF-approximable relation from $(U_{1}, R_{1}, \mathcal{F}_{1})$ to $(U_{2}, R_{2}, \mathcal{F}_{2})$.
\end{dn}

\begin{pn}\label{pn1-CF-apprel} Let $\mathrel{\Theta}$ be a CF-approximable relation from  $(U_{1}, R_{1}, \mathcal{F}_{1})$ to  $(U_{2}, R_{2}, \mathcal{F}_{2})$. Then for all $ F\in \mathcal{F}_{1}, G\in \mathcal{F}_{2}$, the following statements are equivalent:
\begin{enumerate}
\item[$(1)$] $F\mathrel{\Theta} G$;

\item[$(2)$] There exists $F^{\prime}\in \mathcal{F}_{1}$ such that $F^{\prime}\subseteq \overline{R_{1}}(F)$ and $F^{\prime}\mathrel{\Theta} G$;

\item[$(3)$] There exists $G^\prime\in \mathcal{F}_{2}$ such that $F\mathrel{\Theta} G^\prime$ and $G\subseteq\overline{R_{2}}(G^\prime)$;

\item[$(4)$] There exist $F^{\prime}\in \mathcal{F}_{1}$ and  $G^{\prime}\in \mathcal{F}_{2}$ such that $F^{\prime}\subseteq \overline{R_{1}}(F)$, $G\subseteq \overline{R_{2}}(G^{\prime})$ and $F^{\prime}\mathrel{\Theta} G^{\prime}$.
\end{enumerate}
\end{pn}
\begin{proof}
Follows directly from Definition \ref{dn-CF-app}(2)-(4).\end{proof}

Let $\mathrel{\Theta}$ be  a CF-approximable relation from  $(U_{1}, R_{1}, \mathcal{F}_{1})$ to $(U_{2}, R_{2}, \mathcal{F}_{2})$.
 For all $F\in \mathcal{F}_{1}$, set $\widetilde{\mathrel{\Theta}}(F)=\bigcup\{\overline{R_{2}}(G)\mid F\mathrel{\Theta} G\mbox{\ and\ }G\in \mathcal{F}_{2}\}$.  Define a map $f_{\mathrel{\Theta}}: \mathfrak{C}(U_{1}, R_{1}, \mathcal{F}_{1})\longrightarrow \mathcal{P}(U_{2})$ such that for all $ E\in \mathfrak{C}(U_{1}, R_{1}, \mathcal{F}_{1})$, $f_{\mathrel{\Theta}}(E)=\bigcup\{\widetilde{\mathrel{\Theta}}(F)\mid F\subseteq E\mbox{\ and\ }F\in \mathcal{F}_{1}\}$.

\begin{pn}\label{pn2-CF-apprel} Let $\mathrel{\Theta}$ be a CF-approximable relation from  $(U_{1}, R_{1}, \mathcal{F}_{1})$ to   $(U_{2}, R_{2}, \mathcal{F}_{2})$, $F\in \mathcal{F}_{1}$, $E\in \mathfrak{C}(U_{1}, R_{1}, \mathcal{F}_{1})$.  Then the following hold:
\begin{enumerate}
\item[$(1)$] $\{\overline{R_{2}}(G)\mid F\mathrel{\Theta} G\mbox{\ and\ }G\in \mathcal{F}_{2}\}$ is directed;

\item[$(2)$] $\widetilde{\mathrel{\Theta}}(F)\in \mathfrak{C}(U_{2}, R_{2}, \mathcal{F}_{2})$;

\item[$(3)$] For any $F\in \mathcal{F}_{1}$, $f_{\mathrel{\Theta}}(\overline{R_{1}}(F))=\widetilde{\mathrel{\Theta}}(F)$;

\item[$(4)$]   $\{\widetilde{\mathrel{\Theta}}(F)\mid F\subseteq E, F\in \mathcal{F}_{1}\}$ is directed  and $f_{\mathrel{\Theta}}(E)\in \mathfrak{C}(U_{2}, R_{2}, \mathcal{F}_{2})$.
\end{enumerate}
\end{pn}
\begin{proof} (1) By Definition \ref{dn-CF-app}(1), we know that $\{\overline{R_{2}}(G)\mid F\mathrel{\Theta} G, G\in \mathcal{F}_{2}\}$ is not empty. By Definition \ref{dn-CF-app}(5) and
Lemma \ref{lm2-tr-up}, we know that $\{\overline{R_{2}}(G)\mid F\mathrel{\Theta} G, G\in \mathcal{F}_{2}\}$ is directed.

(2) Follows directly from (1) and Proposition \ref{pn2-cf-clo}(1).

(3) It is easy to check that
$$\begin{array}{lll}
f_{\mathrel{\Theta}}(\overline{R_{1}}(F))&=&\bigcup\{\widetilde{\mathrel{\Theta}}(F^{\prime})\mid F^{\prime}\subseteq \overline{R_{1}}(F), F^{\prime}\in \mathcal{F}_{1}\}\\
&=&\bigcup\{\overline{R_{2}}(G)\mid F^{\prime}\in \mathcal{F}_{1}, G\in \mathcal{F}_{2}, F^{\prime}\subseteq \overline{R_{1}}(F)\mbox{~and~}F^{\prime}\mathrel{\Theta} G\}~~(\mbox{by the definition of~} \widetilde{\mathrel{\Theta}}(F^{\prime}))\\
&=&\bigcup\{\overline{R_{2}}(G)\mid G\in \mathcal{F}_{2}, F\mathrel{\Theta} G\}~~(\mbox{by Proposition~}\ref{pn1-CF-apprel}) \\
&=&\widetilde{\mathrel{\Theta}}(F)~~(\mbox{by the definition of~}\widetilde{\mathrel{\Theta}}(F)).
\end{array}$$

(4) Firstly,  we show that $\{\widetilde{\mathrel{\Theta}}(F)\mid F\subseteq E, F\in \mathcal{F}_{1}\}$ is directed. Let $F_{1}, F_{2}\in \mathcal{F}_{1}$. If $F_{1}, F_{2}\subseteq E$, then by Proposition \ref{pn3-cf-cl}(4),  there exists $F_{3}\in \mathcal{F}_{1}$ such that $F_{1}\cup F_{2}\subseteq \overline{R_{1}}(F_{3})\subseteq E$.  And by Definition \ref{dn-CF-app}(2),  it is easy to deduce that $\widetilde{\mathrel{\Theta}}(F_{1})$, $\widetilde{\mathrel{\Theta}}(F_{2})\subseteq \widetilde{\mathrel{\Theta}}(F_{3})$. This shows that the family $\{\widetilde{\mathrel{\Theta}}(F)\mid F\subseteq E, F\in \mathcal{F}_{1}\}$ is directed. Noticing that $f_{\mathrel{\Theta}}(E)=\bigcup\{\widetilde{\mathrel{\Theta}}(F)\mid F\subseteq E, F\in \mathcal{F}_{1}\}$, by (2) and Proposition
\ref{pn2-cf-clo}(3), we have that $f_{\mathrel{\Theta}}(E)\in \mathfrak{C}(U_{2}, R_{2}, \mathcal{F}_{2})$. \end{proof}

Proposition \ref{pn2-CF-apprel} shows that $f_{\mathrel{\Theta}}$ can be seen as a map from $\mathfrak{C}(U_{1}, R_{1}, \mathcal{F}_{1})$ to $\mathfrak{C}(U_{2}, R_{2}, \mathcal{F}_{2})$.

\begin{pn}\label{pn3-CF-apprel} Let $\mathrel{\Theta}$ be  a     CF-approximable relation from        $(U_{1}, R_{1}, \mathcal{F}_{1})$ to    $(U_{2}, R_{2}, \mathcal{F}_{2})$ , $ F\in \mathcal{F}_{1}$, $ G\in \mathcal{F}_{2}$.
 Then $G\subseteq \widetilde{\mathrel{\Theta}}(F)\Leftrightarrow F\mathrel{\Theta} G$.
\end{pn}
\begin{proof} It is easy to check that
$$\begin{array}{lll}
G\subseteq \widetilde{\mathrel{\Theta}}(F)&\Leftrightarrow& G\subseteq \bigcup\{\overline{R_{2}}(G^{\prime})\mid F\mathrel{\Theta} G^{\prime}, G^{\prime}\in \mathcal{F}_{2}\}\\
&\Leftrightarrow&\exists G^{\prime}\in \mathcal{F}_{2}\mbox{~s.t.~}F\mathrel{\Theta} G^{\prime}, G\subseteq\overline{R_{2}}(G^{\prime})~~(\mbox{by Proposition~}
\ref{pn2-CF-apprel}(1)\mbox{~and the finiteness of~}G)\\
&\Leftrightarrow&F\mathrel{\Theta} G~~(\mbox{by Proposition~}
\ref{pn1-CF-apprel}).
\end{array}$$

\vspace{-.3in}\end{proof}
\begin{tm}\label{tm-cfap-sco} Let $\mathrel{\Theta}$ be  a     CF-approximable relation from     $(U_{1}, R_{1}, \mathcal{F}_{1})$ to  $(U_{2}, R_{2}, \mathcal{F}_{2})$.
 Then $f_{\mathrel{\Theta}}$ is a  Scott continuous map from
$\mathfrak{C}(U_{1}, R_{1}, \mathcal{F}_{1})$
 to $\mathfrak{C}(U_{2}, R_{2}, \mathcal{F}_{2})$.
\end{tm}
\begin{proof}  It follows from the definition of $f_{\mathrel{\Theta}}$ that $f_{\mathrel{\Theta}}$ is order preserving.   In order to prove $f_{\mathrel{\Theta}}$ is Scott continuous, by  Proposition \ref{pn2-cf-clo}(3), it suffices to show that for any directed family $\{E_{i}\}_{i\in I}\subseteq \mathfrak{C}(U_{1}, R_{1}, \mathcal{F}_{1})$, we have
$f_{\mathrel{\Theta}}(\bigcup_{i\in I}E_{i})=\bigcup_{i\in I}f_{\mathrel{\Theta}}(E_{i})$. In fact,
$$\begin{array}{lll}
f_{\mathrel{\Theta}}(\bigcup_{i\in I}E_{i})&=&\bigcup\{\widetilde{\mathrel{\Theta}}(F)\mid F\subseteq\bigcup_{i\in I}E_{i}, F\in \mathcal{F}_{1}\}\\
&=&\bigcup_{i\in I}(\bigcup\{\widetilde{\mathrel{\Theta}}(F)\mid F\subseteq E_{i}, F\in \mathcal{F}_{1}\}\\
&&(\mbox{by the finitness of~}F\mbox{ ~and the directedness of ~}\{E_{i}\}_{i\in I})\\
&=&\bigcup_{i\in I}f_{\mathrel{\Theta}}(E_{i})~~(\mbox{by the definition of~}f_{\mathrel{\Theta}}(E_{i})).
\end{array}$$

\vspace{-.45in}
\end{proof}

Theorem \ref{tm-cfap-sco} shows that a CF-approximable relation between  CF-approximation spaces can induce a Scott continuous map between continuous domains. Conversely, a Scott continuous map between  relative continuous domains can also induce a CF-approximable relation between  CF-approximation spaces as follows.

\begin{tm}
 Let $f: \mathfrak{C}(U_{1},  R_{1}, \mathcal{F}_{1})\longrightarrow\mathfrak{C}(U_{2},  R_{2}, \mathcal{F}_{2})$ be a Scott continuous map between CF-approximation spaces $(U_{1}, R_{1}, \mathcal{F}_{1})$ and $(U_{2}, R_{2}, \mathcal{F}_{2})$.
Define $\mathrel{\Theta}_{f}\subseteq\mathcal{F}_{1}\times\mathcal{F}_{2}$
such that\\

$\forall F\in  \mathcal{F}_{1}, G\in \mathcal{F}_{2}, F\mathrel{\Theta}_{f} G\Leftrightarrow G\subseteq f(\overline{R_{1}}(F)).$\\

\noindent Then $\mathrel{\Theta}_{f}$ is a CF-approximable relation from $(U_{1},  R_{1}, \mathcal{F}_{1})$ to $(U_{2},  R_{2}, \mathcal{F}_{2})$.
\end{tm}
\begin{proof}  It follows from  $f(\overline{R_{1}}(F))\in\mathfrak{C}(U_{2},  R_{2}, \mathcal{F}_{2})$ and Definition $\ref{dn-cf-cl}$ that $\mathrel{\Theta}_{f}$ satisfies Definition \ref{dn-CF-app}(1).

To check that $\mathrel{\Theta}_{f}$ satisfies Definition
\ref{dn-CF-app}(2), let  $F, F^{\prime}\in \mathcal{F}_{1}$, $G\in  \mathcal{F}_{2}$. Then
$$\begin{array}{lll}
&&F\subseteq \overline{R_{1}}(F^{\prime}), F\mathrel{\Theta}_{f} G\\
&\Rightarrow&F\subseteq \overline{R_{1}}(F^{\prime}), G\subseteq f(\overline{R_{1}}(F))~~(\mbox{by the definiton of~}\mathrel{\Theta}_{f})\\
&\Rightarrow&G\subseteq f(\overline{R_{1}}(F))\subseteq f(\overline{R_{1}}(F^{\prime}))~~(\mbox{by  Lemma\ }\ref{lm2-tr-up}\mbox{~and  the order preservation of~}f)\\
&\Rightarrow&G\subseteq f(\overline{R_{1}}(F^{\prime}))\Leftrightarrow F^{\prime}\mathrel{\Theta}_{f} G.
\end{array}$$

 To check that $\mathrel{\Theta}_{f}$ satisfies Definition
\ref{dn-CF-app}(3), let
 $F\in \mathcal{F}_{1}$, $G, G^{\prime}\in  \mathcal{F}_{2}$. Then
$$\begin{array}{lll}
&&F\mathrel{\Theta}_{f} G, G^\prime\subseteq\overline{R_{2}}(G)\\
&\Rightarrow&G\subseteq f(\overline{R_{1}}(F)), G^\prime\subseteq\overline{R_{2}}(G)\\
&\Rightarrow& G^\prime\subseteq\overline{R_{2}}(G)\subseteq f(\overline{R_{1}}(F))~~(\mbox{by Proposition\  }\ref{pn2-cf-clo}(2)\mbox{~and~}f(\overline{R_{1}}(F))\in \mathfrak{C}(U_{2},  R_{2}, \mathcal{F}_{2}))\\
&\Rightarrow& G^\prime\subseteq f(\overline{R_{1}}(F))\Leftrightarrow F\mathrel{\Theta}_{f} G^{\prime}.
\end{array}$$

 To check that $\mathrel{\Theta}_{f}$ satisfies Definition
\ref{dn-CF-app}(4), let $F\in \mathcal{F}_{1}$,  $G\in \mathcal{F}_{2}$. If $F\mathrel{\Theta}_{f} G$, then $G\subseteq f(\overline{R_{1}}(F))$. By Proposition \ref{pn2-cf-clo}, \ref{pn3-cf-cl}(2) and the Scott continuity of $f$, we know that

\centerline{
$(\ast)$ \indent \indent $f(\overline{R_{1}}(F))
=\bigcup\{f(\overline{R_{1}}(F^{\prime}))\mid F^{\prime}\subseteq \overline{R_{1}}(F), F^{\prime}\in \mathcal{F}_{1}\}$.}
\noindent Therefore we have that
$$\begin{array}{lll}
G\subseteq f(\overline{R_{1}}(F))&\Rightarrow &\exists G^{\prime}\in \mathcal{F}_{2}, \mbox{~s.t.~}G\subseteq \overline{R_{2}}(G^{\prime}) \mbox{~and~} G^{\prime}\subseteq f(\overline{R_{1}}(F))~~(\mbox{by Definition~}\ref{dn-cf-cl})\\

&\Rightarrow& \exists F^{\prime}\in \mathcal{F}_{1},  G^{\prime}\in \mathcal{F}_{2}, \mbox{~s.t.~}G\subseteq \overline{R_{2}}(G^{\prime}), F^{\prime}\subseteq\overline{R_{1}}(F)\mbox{~and~} G^{\prime}\subseteq f(\overline{R_{1}}(F^{\prime}))\\
&&(\mbox{by equation~}(\ast)\mbox{~and the finiteness of~}G^{\prime})\\
&\Rightarrow& \exists F^{\prime}\in \mathcal{F}_{1},  G^{\prime}\in \mathcal{F}_{2},\mbox{~s.t.~}G\subseteq \overline{R_{2}}(G^{\prime}), F^{\prime}\subseteq\overline{R_{1}}(F)\mbox{~and~} F^{\prime} \mathrel{\Theta}_{f}G^{\prime}\\
&&(\mbox{by the definition of~}\mathrel{\Theta}_{f}).
\end{array}$$

To check that $\mathrel{\Theta}_{f}$ satisfies Definition
\ref{dn-CF-app}(5), let
 $F\in \mathcal{F}_{1}$, $G_{1}, G_{2}\in \mathcal{F}_{2}$. If $F\mathrel{\Theta}_{f} G_{1}$ and $F\mathrel{\Theta}_{f} G_{2}$, then $G_{1}\cup G_{2}\subseteq f(\overline{R_{1}}({F}))$. By Definition \ref{dn-cf-cl} and $f(\overline{R_{1}}({F}))\in \mathfrak{C}(U_{2},  R_{2}, \mathcal{F}_{2})$,  there exists $G_{3}\in \mathcal{F}_{2}$ such that $G_{1}\cup G_{2}\subseteq \overline{R_{2}}(G_{3})\subseteq f(\overline{R_{1}}({F}))$ and $G_{3}\subseteq f(\overline{R_{1}}({F}))$. So $F\mathrel{\Theta}_{f} G_{3}$, showing that $\mathrel{\Theta}_{f}$ satisfies Definition
\ref{dn-CF-app}(5).
\end{proof}

\begin{tm}\label{1-1-apsc}
 Let $f: \mathfrak{C}(U_{1},  R_{1}, \mathcal{F}_{1})\longrightarrow\mathfrak{C}(U_{2},  R_{2}, \mathcal{F}_{2})$ be a Scott continuous map between CF-approximation spaces $(U_{1}, R_{1}, \mathcal{F}_{1})$ and $(U_{2}, R_{2}, \mathcal{F}_{2})$,
 $\mathrel{\Theta}$  a CF-approximable relation from $(U_{1}, R_{1}, \mathcal{F}_{1})$ to  $(U_{2}, R_{2}, \mathcal{F}_{2})$. Then
$\mathrel{\Theta}_{f_{\mathrel{\Theta}}}=\mathrel{\Theta}$ and $f_{\mathrel{\Theta}_{f}}=f$.
\end{tm}
\begin{proof}
Let $F\in \mathcal{F}_{1}, G\in \mathcal{F}_{2}$. Then  by Propositions \ref{pn2-CF-apprel}(3) and \ref{pn3-CF-apprel}, we have
$$(F, G)\in \mathrel{\Theta}_{f_{\mathrel{\Theta}}}\Leftrightarrow  G\subseteq f_{\mathrel{\Theta}}(\overline{R_{1}}({F}))=\widetilde{\mathrel{\Theta}}(F)
\Leftrightarrow (F, G)\in \mathrel{\Theta},$$
 showing that $\mathrel{\Theta}_{f_{\mathrel{\Theta}}}=\mathrel{\Theta}$.

For any $E\in \mathfrak{C}(U_{1},  R_{1}, \mathcal{F}_{1})$, we have that
 $$\begin{array}{lll}
 f_{\mathrel{\Theta}_{f}}(E)&=&\bigcup\{\widetilde{\mathrel{\Theta}_{f}}(F)\mid F\subseteq E\mbox{~and~}F\in \mathcal{F}_{1}\}\\
&=&\bigcup\{\overline{R_{2}}(G)\mid F\in \mathcal{F}_{1}, G\in \mathcal{F}_{2}, F\subseteq E\mbox{~and~}F\mathrel{\Theta}_{f}G\}\\
&=&\bigcup\{\overline{R_{2}}(G)\mid F\in \mathcal{F}_{1}, G\in \mathcal{F}_{2}, F\subseteq E\mbox{~and~}G \subseteq f(\overline{R_{1}}({F}))\}\\
&=&\bigcup\{f(\overline{R_{1}}({F}))\mid F\in \mathcal{F}_{1},  F\subseteq E\}~~(\mbox{by Proposition~}\ref{pn3-cf-cl}(2))\\
&=&f(\bigcup\{\overline{R_{1}}({F})\mid F\in \mathcal{F}_{1},  F\subseteq E\})~~(\mbox{by  Scott continuity of~}f)\\
&=&f(E).
\end{array}$$
This shows that $f_{\mathrel{\Theta}_{f}}=f$. \end{proof}

Given a CF-approximation space $(U, R, \mathcal{F})$,  define the  identity on $(U, R, \mathcal{F})$ to be a binary relation $\operatorname{Id}_{(U, R, \mathcal{F})}\subseteq\mathcal{F}\times\mathcal{F}$ such that for all $F, G\in \mathcal{F}, (F,  G)\in \operatorname{Id}_{(U, R, \mathcal{F})}\Leftrightarrow G\subseteq\overline{R}(F).$


Let $(U_{1},  R_{1}, \mathcal{F}_{1})$, $(U_{2},  R_{2}, \mathcal{F}_{2})$, $(U_{3},  R_{3}, \mathcal{F}_{3})$ be CF-approximation spaces, $\mathrel{\Theta}\subseteq\mathcal{F}_{1}\times\mathcal{F}_{2}$, $\mathrel{\Upsilon}\subseteq\mathcal{F}_{2}\times\mathcal{F}_{3}$ be CF-approximable relations. Define $\mathrel{\Upsilon}\circ\mathrel{\Theta}\subseteq\mathcal{F}_{1}\times\mathcal{F}_{3}$, the composition of $\mathrel{\Upsilon}$ and $\mathrel{\Theta}$ by that
for any $F_{1}\in \mathcal{F}_{1}, F_{3}\in \mathcal{F}_{3}$, $(F_{1}, F_{3})\in \mathrel{\Upsilon}\circ\mathrel{\Theta}$ iff there exists $F_{2}\in \mathcal{F}_{2}$ satisfying $(F_{1}, F_{2})\in \mathrel{\Theta}$  and $(F_{2}, F_{3})\in \mathrel{\Upsilon}$.

It is a routine work to check that $\operatorname{Id}_{(U, R, \mathcal{F})}$ is a CF-approximable relation from $(U, R, \mathcal{F})$ to itself. Thus, CF-approximation spaces as  objects and CF-approximable relations as morphisms with the  identities and  compositions defined above, form a category, and denoted by {\bf CF-GA}. 

Let  {\bf CDOM} be the category of continuous domains and  Scott continuous maps. We next show that categories {\bf CF-GA} and  {\bf CDOM} are equivalent.

\begin{lm}{\em(\cite{Well})}\label{lm-cat-eqv} Let $\mathcal{C}, \mathcal{D}$ be two categories. If there is a functor $\Phi: \mathcal{C}\longrightarrow \mathcal{D}$ such that

$(1)$ $\Phi$ is  full, namely, for all $A, B\in ob(\mathcal{C})$, $g\in Mor_{\mathcal{D}}(\Phi(A), \Phi(B))$,  there is  $f\in Mor_{\mathcal{C}}(A, B)$ such that $\Phi(f)=g$;

$(2)$ $\Phi$ is faithful, namely, for all $A, B\in ob(\mathcal{C})$, $f, g\in Mor_{\mathcal{C}}(A, B)$, if $f\neq g$, then $\Phi(f)\neq \Phi(g)$;

$(3)$ for all $B\in ob(\mathcal{D})$, there is $A\in ob(\mathcal{C})$ such that $\Phi(A)\cong B$,

\noindent then $\mathcal{C}$ and $\mathcal{D}$ are  equivalent.
\end{lm}

\begin{tm}
The categories {\bf CF-GA}  and  {\bf CDOM} are equivalent.
\end{tm}
\begin{proof} 
Define $\Psi:$  {\bf CF-GA}$\to$ {\bf CDOM} such that for all $ (U, R, \mathcal{F})\in ob$({\bf CF-GA}), $\Psi((U, R, \mathcal{F}))=(\mathfrak{C}(U, R, \mathcal{F}), \subseteq)\in  ob({\bf CDOM})$; for all $ \mathrel{\Theta}\in Mor$({\bf CF-GA}),  $\Psi(\mathrel{\Theta})=f_{\mathrel{\Theta}}\in  Mor({\bf CDOM})$.\\

 Give a CF-approximation space $(U, R, \mathcal{F})$,
  for any $E\in \mathfrak{C}(U, R, \mathcal{F})$, we have
$$\begin{array}{lll}
\Psi(\operatorname{Id}_{(U, R, \mathcal{F})})(E)&=&f_{\operatorname{Id}_{(U, R, \mathcal{F})}}(E)\\
&=&\bigcup\{\overline{R}({G})\mid F, G\in \mathcal{F}, F\subseteq E\mbox{~and~} (F, G)\in \operatorname{Id}_{(U, R, \mathcal{F})}\}\\
&=&\bigcup\{\overline{R}({G})\mid F, G\in \mathcal{F}, F\subseteq E\mbox{~and~} G\subseteq \overline{R}({F})\}\\
&=&\bigcup\{\overline{R}({F})\mid  F\in \mathcal{F}, F\subseteq E\}~~(\mbox{by Proposition~}\ref{pn2-cf-clo}(1)\mbox{~and~}\ref{pn3-cf-cl}(2))\\
&=&E ~~(\mbox{by Proposition~}\ref{pn3-cf-cl}(2))\\
&=&id_{\mathfrak{C}(U, R, \mathcal{F})}(E).
\end{array}$$
This shows that $\Psi(\operatorname{Id}_{(U, R, \mathcal{F})})=id_{\mathfrak{C}(U, R, \mathcal{F})}$.\\

Let $\mathrel{\Theta}\subseteq\mathcal{F}_{1}\times\mathcal{F}_{2}$, $\mathrel{\Upsilon}\subseteq\mathcal{F}_{2}\times\mathcal{F}_{3}$ be CF-approximable relations. Then for any $E\in \mathfrak{C}(U, R, \mathcal{F})$, we have
$$\begin{array}{lll}
&&\Psi(\mathrel{\Upsilon})\circ \Psi(\mathrel{\Theta})(E)
=f_{\mathrel{\Upsilon}}(f_{\mathrel{\Theta}}(E))\\
&=&\bigcup\{\widetilde{\mathrel{\Upsilon}}(F)\mid F\subseteq f_{\mathrel{\Theta}}(E), F\in \mathcal{F}_{2}\}\\
&=&\bigcup\{\overline{R_3}({G})\mid F\subseteq f_{\mathrel{\Theta}}(E), F\in \mathcal{F}_{2}, G\in \mathcal{F}_{3}\mbox{~and~}F\mathrel{\Upsilon} G\}~~(\mbox{by the definiton of ~}\widetilde{\mathrel{\Upsilon}})\\
&=&\bigcup\{\overline{R_3}({G})\mid F_{1}\in \mathcal{F}_{1}, F_{1}\subseteq E, G_{1}\in \mathcal{F}_{2}, F_{1}\mathrel{\Theta} G_{1}, F\subseteq \overline{R_2}(G_{1}), F\in
\mathcal{F}_{2}, G\in \mathcal{F}_{3}\mbox{~and~}F\mathrel{\Upsilon} G\}\\
&&(\mbox{by the definiton of\ }f_{\mathrel{\Theta}}(E), \mbox{ Proposition~}\ref{pn2-CF-apprel}(1), \mbox{ Theorem~}\ref{tm-cfap-sco} \mbox{~and  finiteness of members in~}\mathcal{F}_{2})\\
&=&\bigcup\{\overline{R_3}({G})\mid F_{1}\in \mathcal{F}_{1}, F_{1}\subseteq E, G_{1}\in \mathcal{F}_{2}, F_{1}\mathrel{\Theta} G_{1},  G\in \mathcal{F}_{3}\mbox{~and~} G_{1}\mathrel{\Upsilon} G\}\\
&&(\mbox{by~}F\subseteq \overline{R_2}(G_{1}), F\mathrel{\Upsilon} G, \mbox{~and Definition~}\ref{dn-CF-app}(2))\\
&=&\bigcup\{\overline{R_3}({G})\mid F_{1}\in \mathcal{F}_{1}, F_{1}\subseteq E,   G\in \mathcal{F}_{3}\mbox{~and~} (F_{1}, G)\in\mathrel{\Upsilon}\circ\mathrel{\Theta}\}~~(\mbox{by~}F_{1}\mathrel{\Theta} G_{1}\mbox {~and~} G_{1}\mathrel{\Upsilon} G)\\
&=&\bigcup\{\widetilde{\mathrel{\Upsilon}\circ\mathrel{\Theta}}( F_{1})\mid F_{1}\in \mathcal{F}_{1}, F_{1}\subseteq E\}~~(\mbox{by the definition of~}\widetilde{\mathrel{\Upsilon}\circ\mathrel{\Theta}}( F_{1}))\\
&=&f_{\mathrel{\Upsilon}\circ\mathrel{\Theta}}(E)=\Psi(\mathrel{\Upsilon}\circ\mathrel{\Theta})(E).
\end{array}$$
This shows that $\Psi(\mathrel{\Upsilon})\circ \Psi(\mathrel{\Theta})=\Psi(\mathrel{\Upsilon}\circ \mathrel{\Theta})$, and thus $\Psi$ is a functor.

To show that {\bf CF-GA} is equivalent to {\bf CDOM}, it suffices to check that $\Psi$ satisfies the three conditions in Lemma \ref{lm-cat-eqv}.

Let $(U_{1},  R_{1}, \mathcal{F}_{1})$, $(U_{2},  R_{2}, \mathcal{F}_{2})$ be CF-approximation spaces,
 $\mathrel{\Theta}_{1}, \mathrel{\Theta}_{2}$ be  CF-approximable relations from  $(U_{1}, R_{1}, \mathcal{F}_{1})$ to $(U_{2}, R_{2}, \mathcal{F}_{2})$. If $\mathrel{\Theta}_{1}\neq\mathrel{\Theta}_{2}$, then by Theorem \ref{tm-rep} we know that
$\mathrel{\Theta}_{1}=\mathrel{\Theta}_{f_{\mathrel{\Theta}_{1}}}\neq\mathrel{\Theta}_{f_{\mathrel{\Theta}_{2}}}=\mathrel{\Theta}_{2}.$
Thus $f_{\mathrel{\Theta}_{1}}\neq f_{\mathrel{\Theta}_{2}}$, showing that $\Psi$ is faithful.

Let $f: \mathfrak{C}(U_{1},  R_{1}, \mathcal{F}_{1})\longrightarrow \mathfrak{C}(U_{2},  R_{2}, \mathcal{F}_{2})$ be a Scott continuous map, by Theorem \ref{tm-rep}, there is $\mathrel{\Theta}_{f}\in  Mor$({\bf CF-GA}) such that $\Psi(\mathrel{\Theta}_{f})=f_{\mathrel{\Theta}_{f}}=f$, showing that $\Psi$ is full.

It is clear by Theorem \ref{tm-rep} that $\Psi$ satisfies the condition (3) in Lemma \ref{lm-cat-eqv}.
\end{proof}

Similarly, we can also establish   categorical equivalences between the  category of    algebraic domains  with Scott continuous maps as morphisms  and the category of topological CF-approximation spaces with CF-approximable relations as morphisms. We leave the details to the interested readers.\\
\section{Conclusions}

This paper generalizes abstract bases to CF-approximation spaces and generalizes the family of round ideals of an abstract basis to the family of CF-closed sets of a CF-approximation space. Thus  a  representation method of various continuous domains including continuous semilattices, continuous sup-semilattices,  continuous domains with  bottom,  continuous domains with  top, continuous lattices, bc-domains, algebraic domains and arithmetic semilattices in the framework of rough set theory is obtained.  CF-approximable relations between CF-approximation spaces are defined, and   categorical equivalence between categories {\bf CF-GA} of CF-approximation spaces and continuous domains and {\bf CDOM} of continuous domains and Scott continuous maps is established. This work strengthens the links among rough set theory, domain theory and topology, and widens the scope of application of rough set theory and domain theory.


\section*{Acknowledgment}  We would like to thank the referees for their valuable suggestions and comments.


\end{document}